\documentclass[a4paper,11pt]{article}

\usepackage{graphicx}
\usepackage{amssymb,amsbsy,amsmath,amsthm,times,mathrsfs,fullpage}
\usepackage{tikz,wrapfig,multirow}
\usetikzlibrary{patterns,arrows,decorations.pathreplacing}
\usepackage[pdftex]{hyperref}
\usepackage{epstopdf}
\usepackage{mathtools}
\hypersetup{plainpages=True, pdfstartview=FitH, bookmarksopen=true, colorlinks=true,linkcolor=blue,citecolor=blue}

\newcommand*\bell{\ensuremath{\boldsymbol\ell}}

\def\pf{\noindent \emph{Proof.}\ }
\def\qed{{\quad\rule{1mm}{3mm}\,}}

\begin{document}

\pgfdeclarelayer{background}
\pgfdeclarelayer{foreground}
\pgfsetlayers{background,main,foreground}

\newtheorem{thm}{Theorem}
\newtheorem{cor}{Corollary}
\newtheorem{lmm}{Lemma}
\newtheorem{conj}{Conjecture}
\newtheorem{pro}{Proposition}
\newtheorem{Def}{Definition}
\theoremstyle{remark}\newtheorem{Rem}{Remark}

\title{Asymptotic Normality for the Size of Graph Tries built from M-ary Tree Labelings}
\author{Michael Fuchs and Tsan-Cheng Yu\thanks{Both authors were partially supported by grant MOST-109-2115-M-004-003-MY2.}\\
    Department of Mathematical Sciences\\
    National Chengchi University\\
    Taipei 116\\
    Taiwan}

\maketitle

\begin{abstract}
Graph tries are a new and interesting data structure proposed by Jacquet in 2014. They generalize the classical trie data structure which has found many applications in computer science and is one of the most popular data structure on words. For his generalization, Jacquet considered the size (or space requirement) and derived an asymptotic expansion for the mean and the variance when graph tries are built from $n$ independently chosen random labelings of a rooted $M$-ary tree. Moreover, he conjectured a central limit theorem for the (suitably normalized) size as the number of labelings tends to infinity. In this paper, we verify this conjecture with the method of moments.
\end{abstract}

\section{Introduction}

{\it Tries}, proposed by de la Briandais in 1959 and named by Fredkin in 1961, are one of the most important data structures on words. They admit many variants and have found numerous applications; see, e.g., \cite{FuHwZa}, \cite{FuLe} or the recent comprehensive book of Jacquet and Szpankowski \cite{JaSz} for in-depth discussions on their applicability.

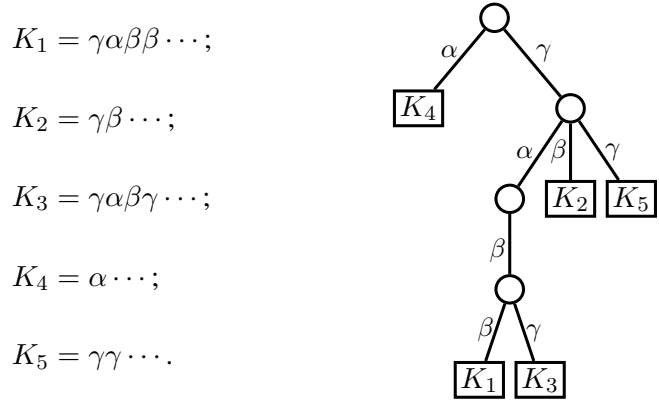
\begin{figure}[h]
\vspace*{0.4cm}
\begin{center}
\begin{tikzpicture}[line width= 0.4mm]
\draw (0cm,0cm) node[circle,draw] (1) {};
\draw (-1cm,-1.2cm) node[inner sep= 2pt,rectangle,draw] (2) {$K_4$};
\draw (1cm,-1.2cm) node[circle,draw] (3) {};
\draw (0.2cm,-2.4cm) node[circle,draw] (4) {};
\draw (1cm,-2.4cm) node[inner sep= 2pt,rectangle,draw] (5) {$K_2$};
\draw (1.8cm,-2.4cm) node[inner sep= 2pt,rectangle,draw] (6) {$K_5$};
\draw (0.2cm,-3.6cm) node[circle,draw] (7) {};
\draw (-0.2cm,-4.8cm) node[inner sep= 2pt,rectangle,draw] (8) {$K_1$};
\draw (0.6cm,-4.8cm) node[inner sep= 2pt,rectangle,draw] (9) {$K_3$};

\draw (1)--(2) node[pos= 0.7,above=0.05cm] {\small{$\alpha$}};
\draw (1)--(3) node[pos= 0.7,above=0.05cm] {\small{$\gamma$}};
\draw (3)--(4) node[pos= 0.85,above=0.1cm] {\small{$\alpha$}};
\draw (3)--(5) node[pos= 0.5,left=-0.1cm] {\small{$\beta$}};
\draw (3)--(6) node[pos= 0.85,above] {\small{$\gamma$}};
\draw (4)--(7) node[pos= 0.6,left=-0.1cm] {\small{$\beta$}};
\draw (7)--(8) node[pos= 1,above=0.12cm] {\small{$\beta$}};
\draw (7)--(9) node[pos= 1,above=0.12cm] {\small{$\gamma$}};

\draw (-6.5cm,-0.3cm) node[right] {$K_1=\gamma\alpha\beta\beta\cdots$;};
\draw (-6.5cm,-1.35cm) node[right] {$K_2=\gamma\beta\cdots$;};
\draw (-6.5cm,-2.4cm) node[right] {$K_3=\gamma\alpha\beta\gamma\cdots$;};
\draw (-6.5cm,-3.45cm) node[right] {$K_4=\alpha\cdots$;};
\draw (-6.5cm,-4.5cm) node[right] {$K_5=\gamma\gamma\cdots$.};
\end{tikzpicture}
\end{center}
\caption{A trie built from five keys with letters from the alphabet ${\mathcal A}=\{\alpha,\beta,\gamma\}$.}\label{trie-fig}
\label{dst-fig}
\end{figure}

We briefly recall the definition of tries which are rooted trees with leaves holding the data that contains infinite strings over a finite alphabet ${\mathcal A}$ of size $A=\vert{\mathcal A}\vert$ as keys. The trie is then an $A$-ary tree which is recursively built from $n$ given keys as follows. For $n=1$, the sole key is placed into the root. For $n\geq 2$, the root is an (empty) internal node and all keys are distributed to the $A$ subtrees of the root according to their first letter; then, the construction of the subtrees proceeds recursively by considering the keys with the first letters removed; see Figure~\ref{trie-fig} for an example.

Note that the above procedure yields the same trie regardless of the order of the $n$ keys. This follows, e.g., from the observation that if one considers the rooted infinite $A$-ary tree with edges labeled by the letters of ${\mathcal A}$ (such that each node has an edge to a child for each letter from ${\mathcal A}$), then a node of this $A$-ary tree is an internal node of the trie if and only if there are at least two keys with prefixes equal to the sequences of letters from the root of the $A$-ary tree to the node.

In order to understand the performance of tries, it is often assumed that the keys are generated by some random process. The simplest of these random processes assumes that the letters of the keys are independent random variables with an identical probability distribution on the alphabet ${\mathcal A}$, i.e., ${\mathbb P}(\alpha)=p_{\alpha}$ for $\alpha\in{\mathcal A}$ where the probabilities $p_{\alpha}$ satisfy $0<p_{\alpha}<1$ and
\[
\sum_{\alpha\in{\mathcal A}}p_{\alpha}=1.
\]
This random model has been used in most of the studies on tries despite it being an oversimplification; see \cite{JaSz}. One reason for this is that it captures many of the phenomena observed in more realistic models, e.g., Markov models (see Leckey et al. \cite{LeNeSz} and references therein).

\vspace*{-0.2cm}
An important shape parameter for tries is the number of internal nodes because it is a measure for the space requirement. If the trie is built from $n$ random keys, we denote the size by $S_n$.

Moments of $S_n$ have been studied in many papers, one of the earliest of which is the paper of Jacquet and R\'{e}gnier \cite{JaRe} (the study of the size over a binary alphabet $\{0,1\}$ with $p_0=p_1=1/2$ goes back even further). Mean and variance have been shown in \cite{JaRe} to be of linear order:
\[
{\mathbb E}(S_n)\sim P_{E}(\log_{a} n)n\qquad\text{and}\qquad{\rm Var}(S_n)\sim P_{V}(\log_{a} n)n,\qquad (n\rightarrow\infty),
\]
where $a>1$ is a suitable constant and $P_{E}(x), P_{V}(x)$ are computable $1$-periodic functions; see Section~\ref{pre} for more details. In addition, also a central limit theorem was proved in \cite{JaRe} (see also \cite{FuLe} and Neininger and R\"{u}schendorf \cite{NeRu}):
\[
\frac{S_n-{\mathbb E}(S_n)}{\sqrt{{\rm Var}(S_n)}}\stackrel{d}{\longrightarrow} N(0,1),
\]
where $\stackrel{d}{\longrightarrow}$ denotes convergence in distribution and $N(0,1)$ is the standard normal distribution.

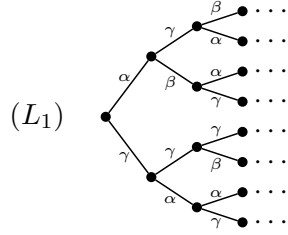
\begin{figure}[h]
\vspace*{0.4cm}
\begin{center}
\begin{tikzpicture}[line width=0.2mm]
\draw (0cm,0cm) node[circle,draw,inner sep=0.04cm,fill] (1) {};
\draw (0.6cm,0.8cm) node[circle,draw,inner sep=0.04cm,fill] (2) {};
\draw (1.2cm,1.2cm) node[circle,draw,inner sep=0.04cm,fill] (3) {};
\draw (1.2cm,0.4cm) node[circle,draw,inner sep=0.04cm,fill] (4) {};
\draw (1.8cm,1.4cm) node[circle,draw,inner sep=0.04cm,fill] (5) {};
\draw (1.8cm,1cm) node[circle,draw,inner sep=0.04cm,fill] (6) {};
\draw (1.8cm,0.6cm) node[circle,draw,inner sep=0.04cm,fill] (7) {};
\draw (1.8cm,0.2cm) node[circle,draw,inner sep=0.04cm,fill] (8) {};
\draw (0.6cm,-0.8cm) node[circle,draw,inner sep=0.04cm,fill] (9) {};
\draw (1.2cm,-0.4cm) node[circle,draw,inner sep=0.04cm,fill] (10) {};
\draw (1.2cm,-1.2cm) node[circle,draw,inner sep=0.04cm,fill] (11) {};
\draw (1.8cm,-0.2cm) node[circle,draw,inner sep=0.04cm,fill] (12) {};
\draw (1.8cm,-0.6cm) node[circle,draw,inner sep=0.04cm,fill] (13) {};
\draw (1.8cm,-1cm) node[circle,draw,inner sep=0.04cm,fill] (14) {};
\draw (1.8cm,-1.4cm) node[circle,draw,inner sep=0.04cm,fill] (15) {};

\draw (2.6cm,1.4cm) node[left] {$\cdots$};
\draw (2.6cm,1cm) node[left] {$\cdots$};
\draw (2.6cm,0.6cm) node[left] {$\cdots$};
\draw (2.6cm,0.2cm) node[left] {$\cdots$};
\draw (2.6cm,-0.2cm) node[left] {$\cdots$};
\draw (2.6cm,-0.6cm) node[left] {$\cdots$};
\draw (2.6cm,-1cm) node[left] {$\cdots$};
\draw (2.6cm,-1.4cm) node[left] {$\cdots$};

\draw (1)--(2) node[pos= 0.4,above] {\tiny{$\alpha$}};
\draw (2)--(3) node[pos= 0.4,above=-0.05cm] {\tiny{$\gamma$}};
\draw (2)--(4) node[pos= 0.4,below=-0.05cm] {\tiny{$\beta$}};
\draw (3)--(5) node[pos= 0.4,above=-0.08cm] {\tiny{$\beta$}};
\draw (3)--(6) node[pos= 0.4,below=-0.08cm] {\tiny{$\alpha$}};
\draw (4)--(7) node[pos= 0.4,above=-0.08cm] {\tiny{$\alpha$}};
\draw (4)--(8) node[pos= 0.4,below=-0.08cm] {\tiny{$\gamma$}};
\draw (1)--(9) node[pos= 0.4,below] {\tiny{$\gamma$}};
\draw (9)--(10) node[pos= 0.4,above=-0.05cm] {\tiny{$\gamma$}};
\draw (9)--(11) node[pos= 0.4,below=-0.05cm] {\tiny{$\alpha$}};
\draw (10)--(12) node[pos= 0.4,above=-0.08cm] {\tiny{$\gamma$}};
\draw (10)--(13) node[pos= 0.4,below=-0.08cm] {\tiny{$\beta$}};
\draw (11)--(14) node[pos= 0.4,above=-0.08cm] {\tiny{$\alpha$}};
\draw (11)--(15) node[pos= 0.4,below=-0.08cm] {\tiny{$\gamma$}};

\draw (-0.9cm,0cm) node {$(L_1)$};
\end{tikzpicture}
\end{center}
\caption{A labeling of a $2$-ary tree with labels from the alphabet ${\mathcal A}=\{\alpha,\beta,\gamma\}$.}\label{M=2}
\end{figure}

Recently, Jacquet \cite{Ja} introduced a fascinating new generalization of tries which he called {\it graph tries} (or G-tries for short). This generalization is built from keys which instead of being strings are now (edge) labelings of a fixed rooted (infinite) graph $G$. (The classical trie is recovered by choosing as $G$ the one-sided infinite path graph.) The concept works for any $G$, however, we will restrict our attention to $M$-ary trees in this paper; see Figure~\ref{M=2} for an example of a labeling of a $2$-ary tree with letters from the alphabet ${\mathcal A}=\{\alpha,\beta,\gamma\}$.

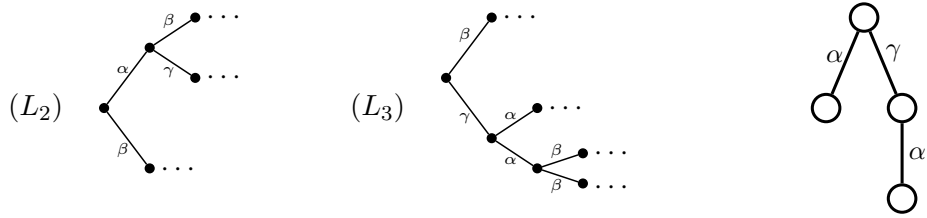
\begin{figure}[h]
\vspace*{0.4cm}
\begin{center}
\begin{tikzpicture}[line width=0.2mm]
\draw (0cm,0cm) node[circle,draw,inner sep=0.04cm,fill] (1) {};
\draw (0.6cm,0.8cm) node[circle,draw,inner sep=0.04cm,fill] (2) {};
\draw (1.2cm,1.2cm) node[circle,draw,inner sep=0.04cm,fill] (3) {};
\draw (1.2cm,0.4cm) node[circle,draw,inner sep=0.04cm,fill] (4) {};
\draw (0.6cm,-0.8cm) node[circle,draw,inner sep=0.04cm,fill] (5) {};

\draw (2cm,0.4cm) node[left] {$\cdots$};
\draw (2cm,1.2cm) node[left] {$\cdots$};
\draw (1.4cm,-0.82cm) node[left] {$\cdots$};

\draw (1)--(2) node[pos= 0.4,above] {\tiny{$\alpha$}};
\draw (2)--(3) node[pos= 0.4,above=-0.05cm] {\tiny{$\beta$}};
\draw (2)--(4) node[pos= 0.4,below=-0.05cm] {\tiny{$\gamma$}};
\draw (1)--(5) node[pos= 0.4,below] {\tiny{$\beta$}};

\draw (-0.9cm,0cm) node {$(L_2)$};

\draw (4.5cm,0.4cm) node[circle,draw,inner sep=0.04cm,fill] (1) {};
\draw (5.1cm,1.2cm) node[circle,draw,inner sep=0.04cm,fill] (2) {};
\draw (5.1cm,-0.4cm) node[circle,draw,inner sep=0.04cm,fill] (3) {};
\draw (5.7cm,0cm) node[circle,draw,inner sep=0.04cm,fill] (4) {};
\draw (5.7cm,-0.8cm) node[circle,draw,inner sep=0.04cm,fill] (5) {};
\draw (6.3cm,-0.6cm) node[circle,draw,inner sep=0.04cm,fill] (6) {};
\draw (6.3cm,-1cm) node[circle,draw,inner sep=0.04cm,fill] (7) {};

\draw (5.9cm,1.2cm) node[left] {$\cdots$};
\draw (6.5cm,0cm) node[left] {$\cdots$};
\draw (7.1cm,-0.6cm) node[left] {$\cdots$};
\draw (7.1cm,-1.05cm) node[left] {$\cdots$};

\draw (1)--(2) node[pos= 0.4,above] {\tiny{$\beta$}};
\draw (1)--(3) node[pos= 0.4,below] {\tiny{$\gamma$}};
\draw (3)--(4) node[pos= 0.4,above=-0.05cm] {\tiny{$\alpha$}};
\draw (3)--(5) node[pos= 0.4,below=-0.05cm] {\tiny{$\alpha$}};
\draw (5)--(6) node[pos= 0.4,above=-0.08cm] {\tiny{$\beta$}};
\draw (5)--(7) node[pos= 0.4,below=-0.08cm] {\tiny{$\beta$}};

\draw (3.6cm,0cm) node {$(L_3)$};

\draw[line width=0.4mm] (10cm,1.2cm) node[circle,draw] (1) {};
\draw[line width=0.4mm] (9.5cm,0cm) node[circle,draw] (2) {};
\draw[line width=0.4mm] (10.5cm,0cm) node[circle,draw] (3) {};
\draw[line width=0.4mm] (10.5cm,-1.2cm) node[circle,draw] (4) {};

\draw[line width=0.4mm] (1)--(2) node[pos=0.9,above=0.2cm] {\small{$\alpha$}};
\draw[line width=0.4mm] (1)--(3) node[pos=0.9,above=0.2cm] {\small{$\gamma$}};
\draw[line width=0.4mm] (3)--(4) node[pos=0.5,right=-0.06cm] {\small{$\alpha$}};
\end{tikzpicture}
\end{center}
\caption{A G-trie (right) built from the three labelings $L_1$ (from Figure~\ref{M=2})$,L_2,$ and $L_3$ of a $2$-ary tree with letters from the alphabet ${\mathcal A}=\{\alpha,\beta,\gamma\}$. For instance, the node of depths $2$ in the G-trie arises from the paths of length $2$ in $L_1$ and $L_3$ which are labeled by $\gamma\alpha$.}\label{G-trie}
\end{figure}

A G-trie is now built from $n$ labelings of the $M$-ary tree as follows. First consider the rooted infinite $M\times A$-ary tree whose nodes correspond to the paths in the $M$-ary tree with all possible labelings, e.g., the root in this $M\times A$-ary tree corresponds to the empty path with an empty labeling and the children of the root correspond to all paths of length $1$ from the root in the $M$-ary tree with all possibilities of assigning labels to these paths. Then, a node in this $M\times A$-ary tree is contained in the G-trie if and only if the path together with its labeling corresponding to the node occurs in at least two of the $n$ given labelings; see Figure~\ref{G-trie} for an example.

As for the classical trie, we assume that all the labels are chosen independently with ${\mathbb P}(\alpha)=p_{\alpha}$ for $\alpha\in{\mathcal A}$ where the probabilities $p_{\alpha}$ are as above. Moreover, we are again interested in the space requirement, i.e., the number of (internal) nodes of a G-trie built from $n$ random labelings which (with a slight abuse of notation) we also denote by $S_n$.

Mean and variance of $S_n$ were investigated in \cite{Ja}. In order to explain the results, we need some notation. First, we order the probabilities of the letters of ${\mathcal A}$ as $p_1,\ldots,p_{A}$. Next, throughout this work, we assume that
\begin{equation}\label{none-exp}
\sum_{j=1}^{A}p_j^2<\frac{1}{M}
\end{equation}
because otherwise ${\mathbb E}(S_n)=\infty$ for $n\geq 2$; see \cite{Ja} where (\ref{none-exp}) was called the {\it non-explosive case}. Finally, let $\rho$ be the  the unique real root with
\[
\sum_{j=1}^{A}p_j^{\rho}=\frac{1}{M}.
\]
Note that $1\leq\rho<2$ because we are in the non-explosive case. Also, $\sigma=1$ only for the classical trie.

Using this notation, it was shown in \cite{Ja} that
\begin{equation}\label{asymp-mean}
{\mathbb E}(S_n)\sim P_{E}(\log_{a} n)n^{\rho},\qquad (n\rightarrow\infty),
\end{equation}
where $a>1$ is a suitable constant and $P_{E}(x)$ is a computable 1-periodic function (depending on $p_1,\ldots,p_A$) which will be given in Section~\ref{pre}. Moreover, for the variance,
\begin{equation}\label{asymp-var}
{\rm Var}(S_n)\sim\begin{cases}P_{V}(\log_{a} n)n^{\rho},&\text{if}\ p_j=1/A\ \text{for}\ 1\leq j\leq A; \\ P_{V}(\log_{a} n)n^{2\rho-1},&\text{otherwise},\end{cases}
\qquad (n\rightarrow\infty),
\end{equation}
where $a$ is as above and $P_{V}(x)$ is another computable 1-periodic function (also depending on $p_1,\ldots,p_A$) which again will be given in Section~\ref{pre}.

The above results are interesting because they show that the size of a G-trie if $\rho>1$ and in the {\it non-uniform case} (i.e., the case where we do not have $p_j=1/A$ for $1\leq j\leq A$) exhibits a very different behavior compared to the size of a classical trie (which is the case $\rho=1$). More precisely, mean and variance are in the above cases of a different order since $2\rho-1>\rho$.

Nevertheless, even in the cases were the variance has a larger order than the mean, the standard deviation is still of a smaller order. Consequently, in \cite{Ja}, a central limit theorem was conjectured. The main result of this paper is the verification of this conjecture.

\begin{thm}\label{main-result}
For the size $S_n$ of a random G-trie built from $n$ labelings of an $M$-ary tree, we have the central limit theorem:
\[
\frac{S_n-{\mathbb E}(S_n)}{\sqrt{{\rm Var}(S_n)}}\stackrel{d}{\longrightarrow} N(0,1).
\]
\end{thm}

This result is interesting because central limit theorems for shape parameters of discrete random structures arising from computer science which have a variance which is significantly larger than the mean are rare; the only other example with a similar behavior for mean and variance can be found in Flajolet et al. \cite{FlSzVa}. Also, in both these examples, most standard methods for proving a central limit theorem do not seem to work (see, e.g.,  Remark~\ref{rem-con-meth} in Section~\ref{non-unif}), leaving only the (more technical) method of moments as a feasible approach; see \cite{FlSzVa} and Section 30 of Billingsley \cite{Bi}.

We conclude the introduction by giving a sketch of the proof of our main result and at the same time outlining the structure of the paper. First, from the definition of random G-tries, $S_n$ satisfies
\begin{equation}\label{dis-rec}
S_n\stackrel{d}{=}\sum_{i=1}^{M}\sum_{j=1}^{A}S_{B_j^{(i)}}^{(i,j)}+1,\qquad (n\geq 2),
\end{equation}
where $S_0=S_1=0$, $S_n^{(i,j)}\stackrel{d}{=}S_n$ for $1\leq i\leq M, 1\leq j\leq A$, and $(B_{1}^{(i)},\ldots,B_{A}^{(i)})$ are independent multinomial distributed random vectors with parameter $(n,p_1,\ldots,p_A)$ for $1\leq i\leq M$. This follows from the fact that $S_n$ is {\it additive} in the sense that it can be computed by first computing it for the subtrees of the root (this gives the terms $S_{B_j^{(i)}}^{(i,j)}$), adding up these contributions, and then increasing this sum by $1$ in order to include the root.

By taking moments on both sides of (\ref{dis-rec}), we see that all (central and non-central) moments satisfy a recurrence of the form
\begin{equation}\label{fund-rec}
a_n=M\sum_{j=1}^{A}\sum_{k=0}^{n}\binom{n}{k}p_j^{k}(1-p_j)^{n-k}a_k+b_n,\qquad (n\geq 2),
\end{equation}
where $a_0=a_1=0$ and $b_n$ with $n\geq 2$ is a given sequence called {\it toll-sequence}.

In \cite{Ja}, the above recurrence for the mean and variance was (asymptotically) solved by using complex-analytic tools such as Mellin transform and analytic depoissonization (for background on these tools see Flajolet et al. \cite{DuFlGo} and Jacquet and Szpankowski \cite{JaSz2}). The results from \cite{Ja} (in a less precise form) were given above; more detailed version will be given in the next section. Moreover, in the next section, we will also give some results concerning the recurrence (\ref{fund-rec}) which exhibits a treshold phenomena: for toll-sequences which have a growth order smaller than $n^{\rho}$, $a_n$ grows like $n^{\rho}$ ({\it small toll-sequence case}); for toll-sequences which have a growth order larger than $n^{\rho}$, $a_n$ grows like $b_n$ ({\it large toll-sequence case}). Here, the large toll-function case will turn out to be crucial for us since for the computation of moments of third and higher order of $S_n$ we will be in this case. In fact, if $\rho>1$ and in the non-uniform case, even the variance falls into the large toll-function case and we will re-derive in this situtation the result from \cite{Ja} in the next section by elementary tools (i.e., tools which are not based on complex analysis).

Then, in Section~\ref{unif} and Section~\ref{non-unif}, we will compute higher moments of $S_n$ again by elementary tools starting from the third moment. The reason why we will consider the uniform case and non-uniform case separately is that these cases behave very differently. In the uniform case, we will see that the computation of the third moment is basically enough to understand the limiting behavior of $S_n$ because with our arguments, the assumptions of the contraction method from \cite{NeRu} can be verified which then immediately implies our desired central limit theorem (this was already suggested in \cite{Ja}). Alternatively, we can also use the method of moments which will be done in Section~\ref{unif} as well. On the other hand, in the non-uniform case, the use of the contraction method is not obvious (see Remark~\ref{rem-con-meth} in Section~\ref{non-unif}) and it seems that only the method of moments works in this case. Technical details for this case are more demanding than in the uniform case and will be given in Section~\ref{non-unif}.

Finally, in Section~\ref{ext}, we will discuss some extensions of our main theorem and in Section~\ref{con} we will give some concluding remarks.

\section{Preliminaries}\label{pre}

This section contains some preliminary results which will be used in the latter sections. Moreover, as mentioned in the introduction, we will recall in more detail the results for the mean and the variance of the size of G-tries from \cite{Ja}. For the sake of clarity, we will divide the section into four paragraphs.

\paragraph{Recurrences for Moments.} Let $\mu_n:={\mathbb E}(S_n)$. Then, by taking expectations on both sides of (\ref{dis-rec}), we obtain that
\[
\mu_n=M\sum_{j=1}^A\sum_{k=0}^n\binom{n}{k}p_j^{k}(1-p_j)^{n-k}\mu_k+1,\qquad (n\geq 2)
\]
with $\mu_0=\mu_1=0$. Note that this is (\ref{fund-rec}) with $b_n=1$ for $n\geq 2$.

Next, we consider central moments. First, note that (\ref{dis-rec}) implies that
\begin{equation}\label{cen-dis-rec}
\left(S_n-\mu_n\right)\stackrel{d}{=}\sum_{i=1}^{M}\sum_{j=1}^{A}\left(S_{B_j^{(i)}}^{(i,j)}-\mu_{B_j^{(i)}}\right)+\Delta_{n,{\mathbf B}^{(1)},\ldots,{\mathbf B}^{(M)}},
\end{equation}
where ${\mathbf B}^{(i)}:=(B_{1}^{(i)},\ldots,B_{A}^{(i)})$ for $1\leq i\leq M$ and
\begin{equation}\label{def-delta}
\Delta_{n,{\bf B}^{(1)},\ldots,{\bf B}^{(M)}}:=1-\mu_n+\sum_{i=1}^{M}\sum_{j=1}^{A}\mu_{B_j^{(i)}}.
\end{equation}
Set
\[
A_{n}^{(m)}:={\mathbb E}\left(S_n-\mu_n\right)^m.
\]
Taking expectations of the $m$-th power on both sides of (\ref{cen-dis-rec}), conditioning on ${\mathbf B}^{(i)}$ on the right-hand side, and expanding gives
\begin{equation}\label{rec-central-mom}
A_n^{(m)}=M\sum_{j=1}^{A}\sum_{k=0}^{n}\binom{n}{k}p_j^{k}(1-p_j)^{n-k}A_j^{(m)}+B_n^{(m)},\qquad (n\geq 2)
\end{equation}
with $A_0^{(m)}=A_1^{(m)}=0$ and
\[
B_{n}^{(m)}:=\sum_{{\mathbf k}^{(1)},\ldots,{\mathbf k}^{(M)}}\left(\prod_{i=1}^{M}\pi_{n,{\mathbf k}^{(i)}}\right)\sum_{\bell^{(1)},\ldots,\bell^{(M)},\ell}\binom{m}{\bell^{(1)},\ldots,\bell^{(M)},\ell}\left(\prod_{i=1}^{M}\prod_{j=1}^{A}A_{k_j^{(i)}}^{(\ell_j^{(i)})}\right)\Delta_{n,{\mathbf k}^{(1)},\ldots,{\mathbf k}^{(M)}}^{\ell},
\]
where the first sum runs over all ${\mathbf k}^{(i)}=(k_1^{(i)},\ldots,k_A^{(i)})\in\{0,\ldots,n\}^{A}$ with $k_1^{(i)}+\cdots+k_A^{(i)}=n$ for $1\leq i\leq M$ and the second sum runs over all $\bell^{(i)}=(\ell_1^{(i)},\ldots,\ell_A^{(i)})\in\{0,\ldots,m-1\}^{A}$ for $1\leq i\leq M$ and $0\leq\ell\leq m$ such that
\begin{equation}\label{sum-cond}
\ell+\sum_{i=1}^{M}\sum_{j=1}^{A}\ell_{j}^{(i)}=m.
\end{equation}
Moreover,
\[
\pi_{n,{\mathbf k}^{(i)}}=\binom{n}{k_1^{(i)},\ldots,k_A^{(i)}}p_1^{k_1^{(i)}}\cdots p_{A}^{k_A^{(i)}}=:\binom{n}{{\bf k}^{(i)}}{\bf p}^{{\bf k}^{(i)}},
\]
where ${\bf p}:=(p_1,\ldots,p_A)$,
\[
\binom{m}{\bell^{(1)},\ldots,\bell^{(M)},\ell}:=\binom{m}{\ell_1^{(1)},\ldots,\ell_A^{(M)},\ell},
\]
and $\Delta_{n,{\mathbf k}^{(1}),\ldots{\mathbf k}^{(M)}}$ is the $\Delta$ from (\ref{def-delta}) but with ${\mathbf B}^{(i)}$ replaced by ${\mathbf k}^{(i)}$ for $1\leq i\leq M$. Note that this shows that the central moments indeed satisfy a recurrence of type (\ref{fund-rec}). (This was claimed at the end of the introduction.)

\paragraph{Asymptotic Expansions for Mean and Variance.} In this paragraph, we will recall more detailed versions of the result for the mean (\ref{asymp-mean}) and variance (\ref{asymp-var}) as given in \cite{Ja}.
In order to state these results, we need the complex roots of the function
\[
P(s):=1-M\sum_{j=1}^{A}p_j^{s}
\]
which satisfy $\Re(s)\leq\rho$. We introduce two notations for the solution set:
\[
{\mathcal S}_{\rho}:=\{s\ :\ P(s)=0\ \text{and}\ \Re(s)=\rho\}
\quad\text{and}\quad
{\mathcal S}_{(\rho-1,\rho]}:=\{s\ :\ P(s)=0\ \text{and}\ \Re(s)\in(\rho-1,\rho]\}.
\]
A great deal about these sets is known; see the deep study of Flajolet et al. \cite{FlRoVa} for a slightly restricted setting which, however, carries over to the solution set of $P(s)$, too.

First, the structure of the set ${\mathcal S}_{\rho}$ depends on a property of the ratios $\log p_i/\log p_j$.
\begin{itemize}
\item[(i)] If the ratios $\log p_i/\log p_j$ are rational for all $1\leq i,j\leq A$, then
\[
{\mathcal S}_{\rho}=\left\{s=\rho+\frac{2k\pi i}{\log a}\ :\ k\in{\mathbb Z}\right\},
\]
where $a>1$ is such that $p_j=a^{-e_j}$ for suitable positive integers $e_j$ for all $1\leq j\leq A$. This is called the {\it periodic case}.
\item[(ii)] If at least one of the ratios $\log p_i/\log p_j$ is irrational for $1\leq i,j\leq A$, then ${\mathcal S}_{\rho}=\{\rho\}$. This is called the {\it aperiodic case}.
\end{itemize}

Secondly, the roots in ${\mathcal S}_{(\rho-1,\rho]}\setminus{\mathcal S}_{\rho}$ behave in more a chaotic way, however, they still exhibit some regularities. For instance, the roots are all simple, are uniformly separated, $1/P(z)$ is bounded provided that $z$ stays uniformly far away from the roots, etc.; see \cite{FlRoVa} for more properties.

Using the above sets, we can now give more detailed versions of (\ref{asymp-mean}) and (\ref{asymp-var}).

\begin{pro}[Jacquet \cite{Ja}]\label{pro-mean}
For the mean of the size $S_n$ of a random G-trie built from $n$ labelings of an $M$-ary tree, we have
\begin{equation}\label{exp-mean}
{\mathbb E}(S_n)=\sum_{\beta\in{\mathcal S}_{(\rho-1,\rho]}}G_{E}(\beta)n^{\beta}+
{\mathcal O}(n^{\rho-1}),
\end{equation}
where the sum and all its term-by-term derivatives are absolute convergent and
\begin{equation}\label{func-mean}
G_E(s)=\frac{(1-s)\Gamma(-s)}{M\sum_{j=1}^{A}p_j^s\log p_j}.
\end{equation}
In particular,
\begin{equation}\label{main-exp-mean}
{\mathbb E}(S_n)\sim\sum_{\beta\in{\mathcal S}_{\rho}}G_E(\beta)n^{\beta}.
\end{equation}
\end{pro}

\begin{Rem}
By setting $a=e$ in the aperiodic case, we recover (\ref{asymp-mean}) with an infinitely divisible periodic function $P_E(x)$.
\end{Rem}

\begin{pro}[Jacquet \cite{Ja}]\label{pro-var}
For the variance of the size $S_n$ of a random G-trie built from $n$ labelings of an $M$-ary tree, we have the following when $\rho>1$.
\begin{itemize}
\item[(i)] In the uniform case, we have
\begin{equation}\label{var-exp-uni}
{\rm Var}(S_n)\sim\sum_{\beta\in{\mathcal S}_{\rho}}G_{V}(\beta)n^{\beta},
\end{equation}
where $G_V(s)$ is a computable function and the sum and all its term-by-term derivates are absolute convergent.
\item[(ii)] In the non-uniform case, we have
\begin{equation}\label{var-exp-non-uni}
{\rm Var}(S_n)\sim\left(\frac{M-1}{MP(2\rho-1)}-1\right)\left(\sum_{\beta\in{\mathcal S_{\rho}}}\beta G_{E}(\beta)n^{\beta-1}\right)^2n,
\end{equation}
where $G_E(x)$ is given in (\ref{func-mean}).
\end{itemize}
\end{pro}

\begin{Rem}
\begin{itemize}
\item[(i)] By again setting $a=e$ in the aperiodic case, we recover (\ref{asymp-var}) with an infinitely divisible periodic function $P_V(x)$.
\item[(ii)] No explicit expression for $G_{V}(s)$ was given in \cite{Ja} where the {\it corrected Poissonized variance approach} from \cite{FuHwZa} was used to prove (\ref{var-exp-uni}). In fact, this approach is capable of giving an explicit (albeit messy) expression for $G_{V}(s)$.
\item[(iii)] A similar result also holds for $\rho=1$ (i.e. for classical tries); see, e.g., \cite{FuHwZa} and \cite{JaRe}.
\end{itemize}
\end{Rem}

In addition to the above results, it was claimed in \cite{Ja} that the mean is of order $n^{\rho}$ and the variance is of order $n^{\rho}$ in the uniform case and of order $n^{2\rho-1}$ in the non-uniform case (for $\rho>1$). Note that this is not obvious from the expressions above. In fact, we have not been able to locate the proofs of these claims in \cite{Ja}. However, we will need these facts (in particular those for the variances) since we will normalize by the right-hand side of (\ref{var-exp-uni}) and (\ref{var-exp-non-uni}) in the proof of our claimed limit law. Therefore, we will give detailed proofs of these claims (amongst other things) in the next two paragraphs.

\paragraph{Small Toll-Sequence Case.} As explained at the end of the introduction, the recurrence (\ref{fund-rec}) exhibits a treshold phenomena. We will give more details in this and the next paragraph. We start with the small toll-sequence case.

First, we point out that it would be possible to prove for the recurrence (\ref{fund-rec}) that $a_n\sim P(\log_{a}n)n^{\rho}$ whenever $b_n$ is small compared to $n^{\rho}$, e.g., $b_n={\mathcal O}(n^{\rho-\epsilon})$ where $\epsilon>0$. (Note that (\ref{asymp-mean}) then immediately follows from this.) Here, $P(x)$ is a 1-periodic function and $a>1$ is as above. However, we do not need such a result in this paper because the small toll-sequence case is just used for the computation of the mean in the general case and the variance for (i) $\rho=1$ and (ii) in the uniform case with $\rho>1$ which was already done in \cite{Ja} with complex-analytic tools that have the additional advantage that they yield precise knowledge of $P(x)$. (Such precise results do not easily follow in the more general situation above.)

Instead, we will recall a result from Schachinger \cite{Scha} (or more precisely, its extension from \cite{FuLe} and Lee \cite{Le}) which shows that if $b_n\geq 0$, then either $a_n$ is constant zero are it grows at least as $n^{\rho}$.

\begin{pro}[Schachinger \cite{Scha}]\label{prop-sch}
Assume that the sequence $a_n$ satisfies (\ref{fund-rec}) with a toll-sequence $b_n$ with (i) $b_n\geq 0$ and (ii) there exist an $n_0\geq 2$ with $b_{n_0}>0$. Then,
\[
a_n=\Omega(n^{\rho}).
\]
\end{pro}

From this, we have the following corollary.
\begin{cor}
The periodic functions $P_E(x)$ for the mean and $P_V(x)$ for the variance in the uniform case and in the non-uniform case with $\rho=1$ are positive.
\end{cor}

\pf The claim about $P_E(x)$ is an immediate consequence of Proposition~\ref{prop-sch}.

As for $P_V(x)$ note that the variance of $S_n$ satisfies (\ref{fund-rec}) with
\begin{equation}\label{toll-var}
b_n=\sum_{{\bf k}^{(1)},\ldots{\bf k}^{(M)}}\left(\prod_{i=1}^{M}\pi_{n,{\bf k}^{(i)}}\right)\Delta_{n,{\bf k}^{(1)},\ldots,{\bf k}^{(M)}}^2;
\end{equation}
see (\ref{rec-central-mom}). This toll-sequence is easily seen to satisfy the assumptions from Proposition~\ref{prop-sch}. Thus,
\[
{\rm Var}(S_n)=\Omega(n^{\rho})
\]
which in the uniform case and in the non-uniform case with $\rho=1$ shows that $P_V(x)>0$.\qed

In fact, the positivity claim for $P_{E}(x)$ can also be proved directly from the explicit expression for $P_E(x)$ which follows from (\ref{func-mean}) and (\ref{main-exp-mean}). This second proof uses ideas from Javanian \cite{Jav} which she introduced to establish a similar claim for a related parameter.

\vspace*{0.2cm}\noindent{\it Second proof that $P_{E}(x)>0$.} First, note that from (\ref{main-exp-mean}), we have
\[
P_E(x)=\sum_{\beta\in{\mathcal S}_{\rho}}G_E(\beta)e^{(\beta-\rho)(\log a)x}.
\]
In the aperiodic case, this sum consists of just one term which is easily seen to be positive.

Thus, we concentrate in the sequel on the periodic case where the above sum becomes
\[
P_E(x)=\sum_{k=-\infty}^{\infty}c_ke^{2\pi k i x}
\]
with $c_k$ given by (\ref{func-mean}) as
\[
c_k=\frac{(1-\rho-\chi_k)\Gamma(-\rho-\chi_k)}{M\sum_{j=1}^{A}p_j^{\rho+\chi_k}\log p_j}=\frac{\Gamma(2-\rho-\chi_k)}{M\left(-\sum_{j=1}^{A}p_j^{\rho}\log p_j\right)(\rho+\chi_k)},
\]
where $\chi_k=2\pi ki/(\log a)$ and we have used that $p_j^{\chi_k}=1$ for $1\leq j\leq A$.

Therefore, it suffices to show that
\[
\sum_{k=-\infty}^{\infty}\frac{\Gamma(2-\rho-\chi_k)}{\rho+\chi_k}e^{2\pi kix}
\]
is positive. In order to establish this claim, observe that
\begin{align*}
\frac{\Gamma(2-\rho-\chi_k)}{\rho+\chi_k}e^{2\pi kix}&=\frac{e^{2\pi kix}}{\rho+\chi_k}\int_{0}^{\infty}e^{-t}t^{1-\rho-\chi_k}{\rm d}t\\
&=e^{2\pi kix}\int_{0}^{\infty}(1-(t+1)e^{-t})t^{-1-\rho-\chi_k}{\rm d}t\\
&=\int_{-\infty}^{\infty}f(e^{u})e^{-\rho u}e^{2\pi kix-u\chi_k}{\rm d}u\\
&=(\log a)\int_{-\infty}^{\infty}f\left(e^{(x-v)(\log a)}\right)e^{-\rho(x-v)(\log a)}e^{2\pi k iv}{\rm d}v,
\end{align*}
where $f(t)=1-(t+1)e^{-t}$ is a positive function on $(0,\infty)$.

Note that the last integral in the computation above is a Fourier transform. Thus, by the Poisson summation formula
\[
\sum_{k=-\infty}^{\infty}\frac{\Gamma(2-\rho-\chi_k)}{\rho+\chi_k}e^{2\pi kix}=(\log a)\sum_{v=-\infty}^{\infty}f\left(e^{(x-v)(\log a)}\right)e^{-\rho(x-v)(\log a)}
\]
which is clearly positive and thus the claim is proved. \qed

\paragraph{Large Toll-Sequence Case.} In this paragraph, we will consider (\ref{fund-rec}) with toll-sequences which grow faster than $n^{\rho}$. This large-toll sequence case will be crucial for the computation of central moments of $S_n$ beyond the second moment and will also enable us to re-prove (\ref{var-exp-non-uni}).

We will show the following asymptotic transfer result.

\begin{pro}\label{asymptotic transfer}
Assume that the sequence $a_n$ satisfies (\ref{fund-rec}) with a toll-sequence $b_n$ with $b_n\sim cn^{\alpha},$ where $\alpha>\rho$ and $c\in{\mathbb R}$. Then,
\[
 a_n\sim\frac{cn^\alpha}{P(\alpha)}.
\]
\end{pro}
\begin{Rem}
\begin{itemize}
\item[(i)] In order to accommodate the case $c=0$, we use (from now on) the convention that $a_n\sim cb_n$ means that $a_n=cb_n+o(b_n)$ (i.e., if $c=0$, the symbol $a_n\sim cb_n$ means that $a_n=o(b_n)$, whereas for $c\ne 0$, the symbol has the usual meaning).
\item[(ii)] This proposition can also be applied term-by-term to $b_n\sim\sum_{\beta}G(\beta)n^{\beta}$ provided that this sum runs over a discrete set of $\beta$'s with $\beta>\rho$ and that it as well as all its term-by-term derivatives are absolute convergent.
\end{itemize}
\end{Rem}

The proof of the above proposition follows by extending the method of proof of Proposition 2 in Hubalek et al. \cite{Hwang}. For the convenience of the reader (and because this proposition will be crucial for the proof of Theorem~\ref{main-result}), we will give a detailed sketch.

First, we need a solution of (\ref{fund-rec}).
\begin{lmm}\label{sumlemma}
Assume that the sequence $a_n$ satisfies (\ref{fund-rec}). Then,
\begin{equation}\label{an-rec}
a_n=\sum_{\ell=0}^{\infty}M^{\ell}\sum_{k_1+\cdots+k_A=\ell}\binom{\ell}{\bf k}\sum_{i=2}^n\binom{n}{i}\left(1-{\bf p}^{{\bf k}}\right)^{n-i}\left({\bf p}^{{\bf k}}\right)^i b_i,
\end{equation}
where ${\bf k}=(k_1,\ldots,k_{A})$, $\binom{\ell}{{\bf k}}$ denotes the multinomial coefficient, and ${\bf p}^{{\bf k}}:=\prod_{j=1}^{A}p_j^{k_j}$.
\end{lmm}
\pf Consider the Poisson-generating functions (as formal power series):
\[
\tilde{f}(z):=e^{-z}\sum_{n=0}^{\infty}\frac{a_n}{n!}z^n\qquad\text{and}\qquad\tilde{g}(z):=e^{-z}\displaystyle\sum_{n=2}^{\infty}\frac{b_n}{n!}z^n.
\]
Then, from (\ref{fund-rec}), we have
\[
\tilde{f}(z) = M\sum_{j=1}^A\tilde{f}(p_jz)+\tilde{g}(z).
\]
Iterating this functional equation gives
\[
\tilde{f}(z)=M^{k+1}\sum_{k_1+\cdots+k_A=k+1}{{k+1}\choose{{\bf k}}}\tilde{f}({\bf p}^{{\bf k}}z)+\sum_{\ell=0}^kM^{\ell}\sum_{k_1+\cdots+k_A=\ell}{{\ell}\choose{{\bf k}}}\tilde{g}({\bf p}^{{\bf k}}z).
\]
By multiplying both sides by $e^{z}$, then taking the $n$-th derivative of both sides and setting $z=0$, we have
\begin{align*}
a_n=M^{k+1}&\sum_{k_1+\cdots+k_A=k+1}{{k+1}\choose{{\bf k}}}\sum_{i=2}^n{{n}\choose{i}}\left(1-{\bf p}^{{\bf k}}\right)^{n-i}\left({\bf p}^{{\bf k}}\right)^ia_i\\
&+\sum_{\ell=0}^kM^{\ell}\sum_{k_1+\cdots+k_A=\ell}{{\ell}\choose{{\bf k}}}\sum_{i=2}^n{{n}\choose{i}}\left(1-{\bf p}^{{\bf k}}\right)^{n-i}\left({\bf p}^{{\bf k}}\right)^ib_i.
\end{align*}
From this the claimed result follows since
\begin{align*}
M^{k+1}&\sum_{k_1+\cdots+k_A=k+1}{{k+1}\choose{{\bf k}}}\sum_{i=2}^n{{n}\choose{i}}\left(1-{\bf p}^{{\bf k}}\right)^{n-i}
\left({\bf p}^{{\bf k}}\right)^ia_i\\
&={\mathcal O}\left(M^{k}\sum_{k_1+\cdots+k_A=k+1}{{k+1}\choose{{\bf k}}}\left({\bf p}^{{\bf k}}\right)^2\right)={\mathcal O}\left(\left(M\sum_{j=1}^{A}p_j^2\right)^k\right)
\end{align*}
approaches 0 as $k\rightarrow\infty$.\qed

\vspace*{0.3cm}
\noindent{\it Proof of Proposition \ref{asymptotic transfer}.} We will split the sum in $k_1, k_2,\ldots, k_A$ in (\ref{an-rec}) into three cases: $n\displaystyle{\bf p}^{{\bf k}}\leq 1$, $1\leq n\displaystyle{\bf p}^{{\bf k}}\leq\log n$, and $n\displaystyle{\bf p}^{{\bf k}}\geq\log n.$

\vspace*{0.2cm}\noindent Case 1: $n{\bf p}^{{\bf k}}\leq 1$. We first estimate the last sum on the right-hand side of (\ref{an-rec}). Uniformly for $\ell\geq 0$ with $k_1+\cdots +k_A=\ell$, we have
\begin{align*}
\sum_{i=2}^{n}{{n}\choose{i}}\left(1-{\bf p}^{{\bf k}}\right)^{n-i}\left({\bf p}^{{\bf k}}\right)^ib_i={\mathcal O}\left(\sum_{i=2}^{n}\frac{i^\alpha}{i!}\left(n{\bf p}^{{\bf k}}\right)^i\right)={\mathcal O}\left(\left(n{\bf p}^{{\bf k}}\right)^2\right)={\mathcal O}\left(\left(n{\bf p}^{{\bf k}}\right)^{\rho + \epsilon}\right),
\end{align*}
for $\epsilon>0$ such that $\rho+\epsilon\leq\min\{2,\alpha\}$. Then,
\begin{align}
\sum_{\ell=0}^{\infty}M^{\ell}\sum_{\substack{k_1+\cdots+k_A=\ell \\ n{\bf p}^{{\bf k}}\leq 1}}&{{\ell}\choose{{\bf k}}}\sum_{i=2}^{n}{{n}\choose{i}}\left(1-{\bf p}^{{\bf k}}\right)^{n-i}\left({\bf p}^{{\bf k}}\right)^ib_i\nonumber\\
=\ &{\mathcal O}\left(\sum_{\ell=0}^{\infty}M^{\ell}\sum_{\substack{k_1+\cdots+k_A=\ell \\ n{\bf p}^{{\bf k}}\leq 1}}{{\ell}\choose{{\bf k}}}\left(n{\bf p}^{{\bf k}}\right)^{\rho+\epsilon}\right).\label{np<1}
\end{align}
Let
\[\hat{p}:= \text{min}\{p_1, \cdots, p_A\}.\]
Then, for any $\ell\geq 0$ with $k_1+\cdots+k_A=\ell$ and $n{\bf p}^{{\bf k}}\leq 1$, we have $\ell\geq -(\log n)/(\log\hat{p})$. Thus,
\begin{align*}
(\ref{np<1})&={\mathcal O}\left(n^{\rho+\epsilon}\sum_{\ell\geq -(\log n)/(\log\hat{p})}\left(M\sum_{j=1}^{A}p_j^{\rho+\epsilon}\right)^{\ell}\right)\\
&={\mathcal O}\left(n^{\rho + \epsilon -\frac{\log \left(M\sum_{j=1}^A p_{j}^{\rho+\epsilon}\right)}{\log \hat{p}}}\right)=o(n^{\alpha}).
\end{align*}

\vspace*{0.1cm}\noindent Case 2: $1\leq n{\bf p}^{{\bf k}}\leq \log n$. Again, we will start with the last sum on the right-hand side of (\ref{an-rec}). Uniformly, for $\ell\geq 1$ with $k_1+\cdots +k_A=\ell$, we have
\begin{equation}
\sum_{i=2}^{n}{{n}\choose{i}}\left(1-{\bf p}^{{\bf k}}\right)^{n-i}\left({\bf p}^{{\bf k}}\right)^ib_i= {\mathcal O}\left((1-{\bf p}^{{\bf k}})^n\sum_{i=2}^{n}\frac{i^\alpha}{i!}\left(\frac{n{\bf p}^{{\bf k}}}{1-{\bf p}^{{\bf k}}}\right)^i\right).\label{1<np<logn}
\end{equation}
Next, by a standard application of the Laplace method or the saddle point method (see Chapter VIII and Appendix B.6 in \cite{FlSe}),
\[
\sum_{j\geq 0}\frac{j^\alpha}{j!}x^{j} = {\mathcal O}(x^\alpha e^x),
\]
for $x\geq 1$. Thus,
\[
(\ref{1<np<logn}) = {\mathcal O}\left((1-{\bf p}^{{\bf k}})^n\left(n{\bf p}^{{\bf k}}\right)^\alpha \exp\left(\frac{n{\bf p}^{{\bf k}}}{1-{\bf p}^{{\bf k}}}\right)\right)={\mathcal O}\left(\left(n{\bf p}^{{\bf k}}\right)^\alpha\right).
\]
The remaining proof proceeds as in Case 1:
\begin{align*}
\sum_{\ell=0}^{\infty}M^{\ell}\sum_{\substack{k_1+\cdots+k_A=\ell \\ 1 \leq n{\bf p}^{{\bf k}} \leq \log n}}&{{\ell}\choose{{\bf k}}}\sum_{i=2}^{n}{{n}\choose{i}}\left(1-{\bf p}^{{\bf k}}\right)^{n-i}\left({\bf p}^{{\bf k}}\right)^ib_i\\
=\ & {\mathcal O}\left(n^{\alpha}\sum_{\ell\geq \hat{\ell}(n)}\left(M\sum_{j=1}^A p_{j}^{\alpha}\right)^{\ell}\right)=o(n^{\alpha}),
\end{align*}
where $\hat{\ell}(n)=-(\log(n/\log n))/(\log\hat{p})$.

\vspace*{0.3cm}\noindent Case 3: $n{\bf p}^{\bf k}\geq\log n$. Also in this case, we will start with the last sum on the right-hand side of (\ref{an-rec}), where now, we will use the local limit theorem of the binomial distribution: for $0<p<1$,
\[
\binom{n}{i}p^{i}(1-p)^{i}=\frac{e^{-t^2/2}}{\sqrt{2\pi np(1-p)}}\left(1+{\mathcal O}\left(\frac{t^3+1}{\sqrt{np(1-p)}}\right)\right)
\]
uniformly for $t=o((p(1-p)n)^{1/6})$, where $i=np+t\sqrt{np(1-p)}$; see Lemma~4 in \cite{Hwang}. Moreover, we will also apply the tail estimates for the binomial distribution from Lemma~5 in \cite{Hwang}. Using these tools, a standard application of the Laplace method gives:
\[
\sum_{i=2}^{n}\binom{n}{i}\left(1-{\bf p}^{\bf k}\right)^{n-i}\left({\bf p}^{\bf k}\right)^{i}b_i\sim c\left(n{\bf p}^{\bf k}\right)^{\alpha}
\]
uniformly for $\ell\geq 0$ with $k_1+\cdots+k_A=\ell$ and $n{\bf p}^{\bf k}\geq\log n$. Thus,
\begin{align*}
\sum_{\ell=0}^{\infty}M^{\ell}\sum_{\substack{k_1+\cdots+k_A=\ell \\ n{\bf p}^{{\bf k}}\geq\log n}}&{{\ell}\choose{{\bf k}}}\sum_{i=2}^{n}{{n}\choose{i}}\left(1-{\bf p}^{{\bf k}}\right)^{n-i}\left({\bf p}^{{\bf k}}\right)^ib_i\\
\sim\ &cn^{\alpha}\sum_{\ell=0}^{\infty}M^{\ell}\sum_{\substack{k_1+\cdots+k_A=\ell \\ n{\bf p}^{{\bf k}}\geq\log n}}{{\ell}\choose{{\bf k}}}\left({\bf p}^{\bf k}\right)^{\alpha}\sim\frac{cn^{\alpha}}{P(\alpha)},
\end{align*}
where in the last step, we used the estimates from Case 1 and Case 2.

Finally, the claimed result follows by combining the three cases.\qed

In the next two sections, we will see that Proposition~\ref{asymptotic transfer} can be used to compute central moments of $S_n$ of order three and higher. In fact, it can be used to compute the variance in the non-uniform case with $\rho>1$ as well which gives an alternative proof of (\ref{var-exp-non-uni}). We will do this next.

First recall that the variance of $S_n$ satisfies (\ref{fund-rec}) with toll-sequence $b_n$ given by (\ref{toll-var}) which can be re-written into
\begin{equation}\label{toll-var-2}
b_n={\mathbb E}\left(\Delta_{n,{\mathbf B}^{(1)},\ldots,{\mathbf B}^{(M)}}\right)^2.
\end{equation}
We will start by deriving the asymptotics of this toll-sequence. Therefore, denote by $E_j^{(i)}$ the event that $\vert B_j^{(i)}-p_{j}n\vert\leq p_{j}n^{2/3}$ and let $E:=\bigcap_{i,j}E_j^{(i)}$. On $E^{c}$, the contribution of the mean on the righ-hand side of (\ref{toll-var-2}) to $b_n$ is exponentially small due to the Chernoff bound for the tail of the binomial distribution. On the other hand, on $E_j^{(i)}$, we obtain by using (\ref{exp-mean}) and Taylor series expansion:
\[
\left(\mu_{B_{j}^{(i)}}\vert E_{j}^{(i)}\right)=\sum_{\beta\in{\mathcal S}_{(\rho-1,\rho]}}G_E(\beta)\left(p_jn\right)^{\beta}+\sum_{\beta\in{\mathcal S_{(\rho-1,\rho]}}}\beta G_E(\beta)\left(p_jn\right)^{\beta-1}\left(B_{j}^{(i)}-p_{j}n\right)+{\mathcal O}(n^{\rho-1}).
\]
Thus, on $E$, we have
\begin{equation}\label{DeltaE}
\left(\Delta_{n,{\mathbf B}^{(1)},\ldots,{\mathbf B}^{(M)}}\vert E\right)=\sum_{\beta\in{\mathcal S}_{(\rho-1,\rho]}}\beta G_E(\beta)\sum_{i=1}^{M}\sum_{j=1}^{A}\left(p_jn\right)^{\beta-1}\left(B_{j}^{(i)}-p_{j}n\right)+{\mathcal O}(n^{\rho-1})
\end{equation}
since
\[
-\sum_{\beta\in{\mathcal S}_{(\rho-1,\rho]}}G_E(\beta)n^{\beta}+\sum_{\beta\in{\mathcal S}_{(\rho-1,\rho]}}G_E(\beta)M\sum_{j=1}^{A}\left(p_jn\right)^{\beta}=-\sum_{\beta\in{\mathcal S}_{(\rho-1,\rho]}}G_E(\beta)n^{\beta}P(\beta)=0.
\]
Consequently, again on $E$,
\[
{\mathbb E}\left(\left(\Delta_{n,{\mathbf B}^{(1)},\ldots,{\mathbf B}^{(M)}}\right)^2\vert E\right)=
M{\mathbb E}\left(\sum_{\beta\in{\mathcal S}_{\rho}}\beta G_E(\beta)\sum_{j=1}^{A}\left(p_jn\right)^{\beta-1}\left(B_{j}^{(1)}-p_{j}n\right)\right)^2+o(n^{2\rho-1}).
\]
Notice that
\begin{align}
{\mathbb E}\Bigg(\sum_{\beta\in{\mathcal S}_{\rho}}\beta G_E(\beta)&\sum_{j=1}^{A}\left(p_jn\right)^{\beta-1}\left(B_{j}^{(1)}-p_{j}n\right)\Bigg)^2\nonumber\\
&={\mathbb E}\left(\sum_{\beta\in{\mathcal S}_{\rho}}\beta G_E(\beta)n^{\beta-1}\sum_{j=1}^{A}p_j^{\rho-1}\left(B_{j}^{(1)}-p_{j}n\right)\right)^2\nonumber\\
&=\left(\sum_{\beta\in{\mathcal S}_{\rho}}\beta G_E(\beta)n^{\beta-1}\right)^2{\mathbb E}\left(\sum_{j=1}^{A}p_j^{\rho-1}\left(B_{j}^{(1)}-p_{j}n\right)\right)^2,\label{simpli}
\end{align}
where we used that $p_j^{\beta}=p_j^{\rho}$ for all $\beta\in{\mathcal S}_{\rho}$ and $1\leq j\leq A$. Next,
\begin{align}
{\mathbb E}\left(\sum_{j=1}^{A}p_j^{\rho-1}\left(B_{j}^{(1)}-p_{j}n\right)\right)^2&=
\sum_{j_1=1}^{A}\sum_{j_2=1}^{A}p_{j_1}^{\rho-1}p_{j_2}^{\rho-1}{\mathbb E}\left(B_{j_1}^{(1)}-p_{j_1}n\right)\left(B_{j_2}^{(1)}-p_{j_2}n\right)\nonumber\\
&=\left(\sum_{j=1}^{A}p_j^{2\rho-1}(1-p_j)-\sum_{j_1\neq j_2}p_{j_1}^{\rho}p_{j_2}^{\rho}\right)n\label{second-form}\\
&=\frac{M-1-MP(2\rho-1)}{M^2}n\label{var-comp}
\end{align}
since
\[
\sum_{j=1}^{A}p_j^{2\rho-1}=\frac{1-P(2\rho-1)}{M}
\]
and
\[
\sum_{j=1}^{A}p_j^{2\rho}+\sum_{j_1\neq j_2}p_{j_1}^{\rho}p_{j_2}^{\rho}=\left(\sum_{j=1}^{A}p_j^{\rho}\right)^2=\frac{1}{M^2}.
\]
Collecting everything and recalling that the contribution to $b_n$ of the expectation on the righ-hand side of (\ref{toll-var-2}) on $E^c$ is exponentially small, we obtain that
\[
b_n={\mathbb E}\left(\Delta_{n,{\mathbf B}^{(1)},\ldots,{\mathbf B}^{(M)}}\right)^2\sim
\frac{M-1-MP(2\rho-1)}{M}\left(\sum_{\beta\in{\mathcal S}_{\rho}}\beta G_E(\beta)n^{\beta-1}\right)^2n.
\]
From this (\ref{var-exp-non-uni}) follows by applying Proposition~\ref{asymptotic transfer}.

The last thing we want to discuss in this section is the positivity of $P_{V}(x)$ in the non-uniform case.
\begin{lmm}
The periodic function $P_V(x)$ for the variance in the non-uniform case with $\rho>1$ is positive.
\end{lmm}

\pf First observe that the coefficient of $n$ in (\ref{second-form}) can be alternatively written as
\[
\sum_{j=1}^{A}p_j^{2\rho-1}(1-p_j)-\sum_{j_1\neq j_2}p_{j_1}^{\rho}p_{j_2}^{\rho}=\sum_{1\leq j_1<j_2\leq A}p_{j_1}p_{j_2}\left(p_{j_1}^{\rho-1}-p_{j_2}^{\rho-1}\right)^2
\]
which shows that this coefficient and consequently also the coefficient of (\ref{var-exp-non-uni}) is positive.

Thus, in order to show the claimed positivity of $P_{V}(x)$, it suffices to show that $Q(x)$ with
\[
\sum_{\beta\in{\mathcal S}_{\rho}}\beta G_{E}(\beta)n^{\beta-1}=Q(\log_a n)n^{\rho-1}
\]
is positive. To do this, we argue as in the second proof that $P_{E}(x)>0$ above.

First, from (\ref{func-mean}) and (\ref{main-exp-mean}), we have
\[
Q(x)=\frac{1}{M\left(-\sum_{j=1}^{A}p_j^{\rho}\log p_j\right)}\sum_{\beta\in{\mathcal S}_{\rho}}\Gamma(2-\beta)e^{(\beta-\rho)(\log a)x}.
\]
In the aperiodic case, the claim is easy; in the periodic case, it is sufficient to prove that
\[
\sum_{k=-\infty}^{\infty}\Gamma(2-\rho-\chi_k)e^{2\pi kix}
\]
is positive, where $\chi_k=2k\pi i/(\log a)$. This follows from
\[
\Gamma(2-\rho-\chi_k)e^{2\pi kix}=(\log a)\int_{-\infty}^{\infty}e^{-e^{(x-v)(\log a)}}e^{-\rho(x-v)(\log a)}e^{2\pi kiv}{\rm d}v,
\]
which follows by a similar computation as in the second proof of the positive of $P_E(x)$, and another application of the Poisson summation formula. \qed

\section{Size: Uniform Case}\label{unif}

We assume throughout this section that $\rho>1$ and that we are in the uniform case even though the arguments below also work for $\rho=1$ (both in the uniform and non-uniform case). However, the latter situation is already covered by previous work; see the introduction.

We first note that in the uniform case, the solution set of $P(s)=1-MA^{1-s}$ becomes much easier. Clearly, $\rho=1+\log_{A}M$ and thus
\[
{\mathcal S}_{\rho}=\left\{s=1+\log_{A}M+\frac{2k\pi i}{\log M}\ :\ k\in{\mathbb Z}\right\}.
\]
Also note that there are no zeros of $P(s)$ with $\Re(s)<\rho$. As a consequence, (\ref{DeltaE}) becomes
\begin{align}
\left(\Delta_{n,{\bf B}^{(1)},\ldots,{\bf B}^{(M)}}\vert E\right)&=\left(\sum_{\beta\in{\mathcal S}_{\rho}}\beta G_{E}(\beta)n^{\beta-1}\right)\left(\sum_{i=1}^M\sum_{j=1}^AA^{1-\rho}\left(B_{j}^{(i)}-\frac{n}{A}\right)\right)+{\mathcal O}(n^{\rho-1})\nonumber\\
&={\mathcal O}(n^{\rho-1})=o(n^{\rho/2})\label{est-delta}
\end{align}
since
\[
A^{1-\rho}\sum_{i=1}^M\sum_{j=1}^A\left(B_j^{(i)}-\frac{n}{A}\right)=A^{1-\rho}\left(Mn-Mn\right)=0
\]
and $\rho-1<\rho/2$. (Because we are in the non-explosive case.)

We now consider the third moment which satisfies (\ref{rec-central-mom}) with
\[
B_{n}^{(3)}=3MA{\mathbb E}\left(A_{B_1^{(1)}}^{(2)}\Delta_{n,{\bf B}^{(1)},\ldots,{\bf B}^{(M)}}\right)+{\mathbb E}\left(\Delta_{n,{\bf B}^{(1)},\ldots,{\bf B}^{(M)}}\right)^3.
\]
From (\ref{est-delta}) and the result for the variance (Proposition~\ref{pro-var}), we have
\[
{\mathbb E}\left(A_{B_1^{(1)}}^{(2)}\Delta_{n,{\bf B}^{(1)},\ldots,{\bf B}^{(M)}}\right)=o(n^{3\rho/2}),
\]
where we used that the contribution of the mean on $E^{c}$ is exponentially small. Likewise,
\[
{\mathbb E}\left(\Delta_{n,{\bf B}^{(1)},\ldots,{\bf B}^{(M)}}\right)^3=o(n^{3\rho/2}).
\]
Thus, $B_n^{(3)}=o(n^{3\rho/2})$ and by Proposition~\ref{asymptotic transfer},
\[
A_n^{(3)}=o(n^{3\rho/2}).
\]

From the above simple argument, we also get a proof of Theorem~\ref{main-result} for the uniform case via the contraction method.

\vspace*{0.3cm}\noindent{\it Proof of Theorem~\ref{main-result} in the uniform case.} We apply Corollary 5.2 in \cite{NeRu} to (\ref{dis-rec}) which implies the claimed result upon checking the conditions of this corollary the most important of which is:
\[
{\mathbb E}\vert\Delta_{n,{\bf B}^{(1)},\ldots,{\bf B}^{(M)}}\vert^3=o(n^{3\rho/2}).
\]
This condition immediately follows from (\ref{est-delta}).\qed

Alternatively, we can compute higher moments which also leads to a proof of Theorem~\ref{main-result} by applying the Fr\'{e}chet-Shohat theorem to the following proposition.

\begin{pro}\label{pro-unif-case}
For $m\geq 0$,
\[
A_n^{(m)}\sim g_m\left(\sum_{\beta\in{\mathcal S}_{\rho}}G_V(\beta)n^{\beta}\right)^{m/2},
\]
where $g_m$ denotes the $m$-th moment of the standard normal distribution, i.e.,
\[
g_m:=\begin{cases} m!/(2^{m/2}(m/2)!), &\text{if}\ m\ \text{is even}; \\ 0, &\text{if}\ m\ \text{is odd}.\end{cases}
\]
\end{pro}
\pf We use induction on $m$. First note that $m=0$ and $m=1$ are trivial, $m=2$ is contained in Proposition~\ref{pro-var}-(i), and $m=3$ was proved above. Next, assume that the statement is true for all $m'<m$. We are going to proof it for $m$.

First, consider $B_n^{(m)}$ which is given by
\[
B_{n}^{(m)}=\sum_{\bell^{(1)},\ldots,\bell^{(M)},\ell}\binom{m}{\bell^{(1)},\ldots,\bell^{(M)},\ell}{\mathbb E}\left(\left(\prod_{i=1}^{M}\prod_{j=1}^{A}A_{B_j^{(i)}}^{(\ell_j^{(i)})}\right)\Delta_{n,{\mathbf B}^{(1)},\ldots,{\mathbf B}^{(M)}}^{\ell}\right).
\]
Note that by (\ref{est-delta}) and the induction hypothesis, whenever $\ell\geq 1$,
\[
{\mathbb E}\left(\left(\prod_{i=1}^{M}\prod_{j=1}^{A}A_{B_j^{(i)}}^{(\ell_j^{(i)})}\right)\Delta_{n,{\mathbf B}^{(1)},\ldots,{\mathbf B}^{(M)}}^{\ell}\Big\vert E\right)=o\left(n^{(\ell+\sum_{i=1}^{M}\sum_{j=1}^{A}\ell_{j}^{(i)})/2}\right)=o(n^{m/2}),
\]
where we used (\ref{sum-cond}). Thus, from the exponential smallness of the above expectation on $E^c$, we have
\[
{\mathbb E}\left(\left(\prod_{i=1}^{M}\prod_{j=1}^{A}A_{B_j^{(i)}}^{(\ell_j^{(i)})}\right)\Delta_{n,{\mathbf B}^{(1)},\ldots,{\mathbf B}^{(M)}}^{\ell}\right)=o(n^{m/2})
\]
which means (in view of the claimed result) that the main contribution to $B_{n}^{(m)}$ comes from the terms with $\ell=0$. For these terms, again by the induction hypothesis and (\ref{sum-cond}) together with Taylor series expansion:
\[
{\mathbb E}\left(\left(\prod_{i=1}^{M}\prod_{j=1}^{A}A_{B_j^{(i)}}^{(\ell_j^{(i)})}\right)\Big\vert E\right)\sim\left(\sum_{\beta\in{\mathcal S}_{\rho}}G_V(\beta)\left(\frac{n}{A}\right)^{\beta}\right)^{m/2}\left(\prod_{i=1}^{M}\prod_{j=1}^Ag_{\ell_j^{(i)}}\right)
\]
and this asymptotic again holds with the conditioning on $E$ removed. Consequently,
\begin{align*}
\sum_{\bell^{(1)},\ldots,\bell^{(M)}}&\binom{m}{\bell^{(1)},\ldots,\bell^{(M)}}{\mathbb E}\left(\left(\prod_{i=1}^{M}\prod_{j=1}^{A}A_{B_j^{(i)}}^{(\ell_j^{(i)})}\right)\right)\\
&\sim\left(\sum_{\beta\in{\mathcal S}_{\rho}}G_V(\beta)n^{\beta}\right)^{m/2}A^{-\rho m/2}\sum_{\bell^{(1)},\ldots,\bell^{(M)}}\binom{m}{\bell^{(1)},\ldots,\bell^{(M)}}\left(\prod_{i=1}^{M}\prod_{j=1}^Ag_{\ell_j^{(i)}}\right)\\
&\sim g_mP(\rho m/2)\left(\sum_{\beta\in{\mathcal S}_{\rho}}G_V(\beta)n^{\beta}\right)^{m/2},
\end{align*}
where in the first step we used that $A^{-\beta}=A^{-\rho}$ for $\beta\in{\mathcal S}_{\rho}$ and in the second step:
\[
\sum_{\bell^{(1)},\ldots,\bell^{(M)}}\binom{m}{\bell^{(1)},\ldots,\bell^{(M)}}\left(\prod_{i=1}^{M}\prod_{j=1}^Ag_{\ell_j^{(i)}}\right)=g_m((MA)^{m/2}-MA)=g_m A^{-\rho m/2}P(\rho m/2)
\]
which is trivial if $m$ is odd (since all terms are zero) and follows easily from the multinomial theorem if $m$ is even. Moreover, we used that $MA=A^{\rho}$ in the last step.

Finally, by applying Proposition~\ref{asymptotic transfer}, we obtain that
\[
A_n^{(m)}\sim g_m\left(\sum_{\beta\in{\mathcal S}_{\rho}}G_V(\beta)n^{\beta}\right)^{m/2}
\]
which proves the claim.\qed

\section{Size: Non-Uniform Case}\label{non-unif}

We assume throughout this section that $\rho>1$ and that we are in the non-uniform case.

The proof of Theorem~\ref{main-result} follows from the the following proposition and the Fr\'{e}chet-Shohat Theorem.

\begin{pro}
For $m\geq 0$,
\[
A_n^{(m)}\sim g_{m}c^{m/2}\left(\sum_{\beta\in{\mathcal S}_{\rho}}\beta G_E(\beta)n^{\beta-1}\right)^mn^{m/2},
\]
where
\[
c:=\frac{M-1}{MP(2\rho-1)}-1
\]
and $g_m$ denotes the $m$-th moment of the standard normal distribution.
\end{pro}
\pf The proof is by induction of $m$. The claim is trivial for $m=0$ and $m=1$ and was proved for $m=2$ in Section~\ref{pre}. Therefore, we assume that the claim is true for all $m'<m$. We are going to establish it for $m$.

First, consider $B_{n}^{(m)}$ which as in proof of Proposition~\ref{pro-unif-case} is given by
\[
B_{n}^{(m)}=\sum_{\bell^{(1)},\ldots,\bell^{(M)},\ell}\binom{m}{\bell^{(1)},\ldots,\bell^{(M)},\ell}{\mathbb E}\left(\left(\prod_{i=1}^{M}\prod_{j=1}^{A}A_{B_j^{(i)}}^{(\ell_j^{(i)})}\right)\Delta_{n,{\mathbf B}^{(1)},\ldots,{\mathbf B}^{(M)}}^{\ell}\right).
\]
Recall the definition of the events $E_j^{(i)}$ in the derivation of the variance from Section~\ref{pre}. On $E_j^{(i)}$, by induction hypothesis and Taylor series expansion, we have
\begin{align*}
\left(A_{B_j^{(i)}}^{(\ell_j^{(i)})}\vert E_j^{(i)}\right)&\sim
g_{\ell_j^{(i)}}c^{\ell_j^{(i)}/2}\left(\sum_{\beta\in{\mathcal S}_{\rho}}\beta G_E(\beta)\left(p_jn\right)^{\beta-1}\right)^{\ell_j^{(i)}}(p_j n)^{\ell_j^{(i)}/2}\\
&\sim
g_{\ell_j^{(i)}}p_j^{(\rho-1/2)\ell_j^{(i)}}c^{\ell_j^{(i)}/2}\left(\sum_{\beta\in{\mathcal S}_{\rho}}\beta G_E(\beta)n^{\beta-1}\right)^{\ell_j^{(i)}}n^{\ell_j^{(i)}/2},
\end{align*}
where we used that $p_j^{\beta}=p_j^{\rho}$ for $\beta\in\mathcal{S}_{\rho}$ and $1\leq j\leq A$. Consequently,
\begin{align}
{\mathbb E}&\Bigg(\Bigg(\prod_{i=1}^{M}\prod_{j=1}^{A}A_{B_j^{(i)}}^{(\ell_j^{(i)})}\Bigg)\Delta_{n,{\mathbf B}^{(1)},\ldots,{\mathbf B}^{(M)}}^{\ell}\Big\vert E\Bigg)\nonumber\\
&\sim\left(\prod_{i=1}^{M}\prod_{j=1}^{A}g_{\ell_j^{(i)}}p_j^{(\rho-1/2)\ell_{j}^{(i)}}\right)c^{(m-\ell)/2}\left(\sum_{\beta\in{\mathcal S}_{\rho}}\beta G_E(\beta)n^{\beta-1}\right)^{m-\ell}n^{(m-\ell)/2}{\mathbb E}\left(\Delta_{n,{\mathbf B}^{(1)},\ldots,{\mathbf B}^{(M)}}^{\ell}\vert E\right),\label{non-uniform-step-1}
\end{align}
where we used (\ref{sum-cond}). Next, from (\ref{DeltaE}),
\begin{equation}\label{non-uniform-step-2}
{\mathbb E}\left(\Delta_{n,{\mathbf B}^{(1)},\ldots,{\mathbf B}^{(M)}}^{\ell}\vert E\right)\sim\left(\sum_{\beta\in{\mathcal S}_{\rho}}\beta G_E(\beta)n^{\beta-1}\right)^{\ell}{\mathbb E}\left(\sum_{i=1}^M\sum_{j=1}^Ap_j^{\rho-1}\left(B_j^{(i)}-p_jn\right)\right)^{\ell},
\end{equation}
where we used a similar computation as in (\ref{simpli}). By expanding the last expectation, we obtain that
\[
{\mathbb E}\left(\sum_{i=1}^M\sum_{j=1}^Ap_j^{\rho-1}\left(B_j^{(i)}-p_jn\right)\right)^{\ell}=\sum_{{\bf k}}\binom{\ell}{{\bf k}}\prod_{i=1}^{M}{\mathbb E}\left(\sum_{j=1}^{A}p_j^{\rho-1}\left(B_j^{(1)}-p_j n\right)\right)^{k_i},
\]
where ${\bf k}=(k_1,\ldots,k_M)$ such that $\sum_{i=1}^{M}k_i=\ell$. Now, from the central limit theorem for the multinomial distribution, we have in distribution and with convergence of all moments:
\[
\frac{\sum_{j=1}^{A}p_j^{\rho-1}\left(B_j^{(1)}-p_j n\right)}{\sqrt{n}}\stackrel{d}{\longrightarrow} N\left(0,\frac{P(2\rho-1)c}{M}\right),
\]
where the variance was computed in (\ref{var-comp}). Thus,
\[
{\mathbb E}\left(\sum_{i=1}^M\sum_{j=1}^Ap_j^{\rho-1}\left(B_j^{(i)}-p_jn\right)\right)^{\ell}\sim\left(\frac{P(2\rho-1)c}{M}\right)^{\ell/2}n^{\ell/2}\sum_{{\bf k}}\binom{\ell}{{\bf k}}\left(\prod_{i=1}^{M}g_{k_i}\right).
\]
Note that
\[
\sum_{{\bf k}}\binom{\ell}{{\bf k}}\left(\prod_{i=1}^{M}g_{k_i}\right)=g_{\ell}M^{\ell/2}
\]
which is trivial if $\ell$ is odd and follows from the multinomial theorem if $\ell$ is even. Consequently,
\[
{\mathbb E}\left(\sum_{i=1}^M\sum_{j=1}^Ap_j^{\rho-1}\left(B_j^{(i)}-p_jn\right)\right)^{\ell}\sim g_{\ell}\left(P(2\rho-1)c\right)^{\ell/2}n^{\ell/2}.
\]
Now, plugging this into (\ref{non-uniform-step-2}) and (\ref{non-uniform-step-2}) in turn into (\ref{non-uniform-step-1}), we obtain that
\begin{align*}
{\mathbb E}\Bigg(\Bigg(\prod_{i=1}^{M}\prod_{j=1}^{A}&A_{B_j^{(i)}}^{(\ell_j^{(i)})}\Bigg)\Delta_{n,{\mathbf B}^{(1)},\ldots,{\mathbf B}^{(M)}}^{\ell}\Big\vert E\Bigg)\\
&\sim\left(\prod_{i=1}^{M}\prod_{j=1}^{A}g_{\ell_j^{(i)}}p_j^{(\rho-1/2)\ell_{j}^{(i)}}\right)g_{\ell}P(2\rho-1)^{\ell/2}c^{m/2}\left(\sum_{\beta\in{\mathcal S}_{\rho}}\beta G_E(\beta)n^{\beta-1}\right)^{m}n^{m/2}
\end{align*}
and the same holds without the conditioning on $E$ since the conditional expectation conditioned on the complement of $E$ is exponentially small. Plugging this into the expression for $B_n^{(m)}$ gives
\[
B_{n}^{(m)}\sim d c^{m/2}\left(\sum_{\beta\in{\mathcal S}_{\rho}}\beta G_E(\beta)n^{\beta-1}\right)^{m}n^{m/2},
\]
where
\[
d=\sum_{\bell^{(1)},\ldots,\bell^{(M)},\ell}g_{\ell}\binom{m}{\bell^{(1)},\ldots,\bell^{(M)},\ell}\left(\prod_{i=1}^{M}\prod_{j=1}^{A}g_{\ell_j^{(i)}}p_j^{(\rho-1/2)\ell_{j}^{(i)}}\right)P(2\rho-1)^{\ell/2}
=g_mP((\rho-1/2)m).
\]
which is trivial if $m$ is odd and follows from the multinomial theorem if $m$ is even. Overall,
\[
B_{n}^{(m)}\sim g_mP((\rho-1/2)m)c^{m/2}\left(\sum_{\beta\in{\mathcal S}_{\rho}}\beta G_E(\beta)n^{\beta-1}\right)^{m}n^{m/2}.
\]

Finally, applying Proposition~\ref{asymptotic transfer} gives
\[
A_n^{(m)}\sim g_mc^{m/2}\left(\sum_{\beta\in{\mathcal S}_{\rho}}\beta G_E(\beta)n^{\beta-1}\right)^{m}n^{m/2}.
\]
This concludes the proof of the induction step.\qed

\begin{Rem}\label{rem-con-meth}
From the above proof, we see that
\[
{\mathbb E}\left(\Delta_{n,{\bf B}^{(1)},\ldots,{\bf B}^{(M)}}^{3}\right)\sim g_3(P(2\rho-1)c)^{3/2}\left(\sum_{\beta\in{\mathcal S}_{\rho}}\beta G_E(\beta)n^{\beta-1}\right)^3n^{3/2}=o\left(n^{3(\rho-1/2)}\right)
\]
since $g_3=0$. In order to be able to apply the contraction method in the non-uniform case, we would need
\[
{\mathbb E}\vert\Delta_{n,{\bf B}^{(1)},\ldots,{\bf B}^{(M)}}\vert^{3}=o\left(n^{3(\rho-1/2)}\right),
\]
however, it is not clear how to prove this (provided that this even holds).
\end{Rem}

\section{Extensions}\label{ext}

The method we used to prove Theorem~\ref{main-result} has some generality and can be applied to other additive parameters of G-tries as well. In this section, we will briefly discuss two such classes of parameters.

\paragraph{External Nodes containing Keys and Empty External Nodes.} The G-trie built from $n$ random keys can be completed by adding external nodes so that all internal nodes of the G-trie have outdegree $MA$. Some of these external nodes will contain data and some not where an external node is of the first type if and only if the path with its labeling represented by the node can be found in exactly one of the $n$ keys.

Denote by $K_n$ and $R_n$ the external nodes of type 1 and type 2. Note that $K_n=n$ if $\rho=1$ and this is the only case where one of these random variables is deterministic. Clearly, $K_n$ and $R_n$ both satisfy the recurrence (\ref{dis-rec}) with the only difference that the initial conditions are different: $K_0=0$ and $K_1=1$ for $K_n$ and $R_0=1$ and $R_1=0$ for $R_n$. In fact, we can more generally define
\[
N_n:=\alpha R_n+\beta K_n+\gamma S_n,
\]
where $\alpha,\beta,\gamma\geq 0$ with at least one of them positive. Then,
\[
N_n\stackrel{d}{=}\sum_{i=1}^{M}\sum_{j=1}^{A}N_{B_j^{(i)}}^{(i,j)}+\gamma,\qquad (n\geq 2),
\]
where notation is as in the introduction and $N_0=\alpha$ and $N_1=\beta$.

Using the tools from \cite{Ja} and from this paper, the following result holds for $N_n$.
\begin{thm}
We have,
\[
{\mathbb E}(N_n)\sim Q_{E}(\log_{a} n)n^{\rho},\qquad (n\rightarrow\infty),
\]
and
\[
{\rm Var}(N_n)\sim\begin{cases}Q_{V}(\log_{a} n)n^{\rho},&\text{if}\ p_j=1/A\ \text{for}\ 1\leq j\leq A; \\ Q_{V}(\log_{a} n)n^{2\rho-1},&\text{otherwise},\end{cases}
\qquad (n\rightarrow\infty),
\]
where $a>1$ and $Q_E(x),Q_V(x)$ are computable 1-periodic function. Moreover,
\[
\frac{N_n-{\mathbb E}(N_n)}{\sqrt{{\rm Var}(N_n)}}\stackrel{d}{\longrightarrow}N(0,1).
\]
\end{thm}

\paragraph{Leaves and Patterns.} Another class of examples arises from counting patterns, e.g., the number of leaves in a G-trie. We use this number as guiding example, although more general patterns can be considered as well.

Denote by $L_n$ the number of leaves of a G-trie built from $n$ random keys, e.g., $L_n=2$ for the G-trie in Figure~\ref{G-trie}. Clearly, $L_n$ is additive since it satisfies the recurrence:
\[
L_n\stackrel{d}{=}\sum_{i=1}^{M}\sum_{j=1}^{A}L_{B_j^{(i)}}^{(i,j)}+T_n,\qquad (n\geq 2),
\]
where $L_0=0$ and $L_1=1$ and
\[
T_n=\begin{cases} 1,&\text{if}\ B_j^{(i)}\in\{0,1\}\ \text{for all}\ 1\leq i\leq M\ \text{and}\ 1\leq j\leq A;\\0,&\text{otherwise}.\end{cases}
\]
Again the tools from \cite{Ja} and this paper can be applied to obtain the following result.
\begin{thm}
We have,
\[
{\mathbb E}(L_n)\sim H_{E}(\log_{a} n)n^{\rho},\qquad (n\rightarrow\infty),
\]
and
\[
{\rm Var}(L_n)\sim\begin{cases}H_{V}(\log_{a} n)n^{\rho},&\text{if}\ p_j=1/A\ \text{for}\ 1\leq j\leq A; \\ H_{V}(\log_{a} n)n^{2\rho-1},&\text{otherwise},\end{cases}
\qquad (n\rightarrow\infty),
\]
where $a>1$ and $H_E(x),H_V(x)$ are computable 1-periodic function. Moreover,
\[
\frac{L_n-{\mathbb E}(L_n)}{\sqrt{{\rm Var}(L_n)}}\stackrel{d}{\longrightarrow}N(0,1).
\]
\end{thm}

\section{Conclusion}\label{con}

In this paper we considered G-tries which are interesting new data structures that were recently proposed by Jacquet who computed the mean and variance of the size of a G-trie. Moreover, he conjectured a central limit theorem (see the conclusion of \cite{Ja}) which we proved in this paper by applying the method of moments. This shape parameter is interesting because its mean and variance in the non-classical and non-uniform case behave very different from the classical case: it is one of the few known examples of a shape parameters of discrete random structures arising from computer science with a variance of a considerable larger asymptotic order than the mean yet still a central limit theorem holds. The proof of the central limit theorem for this parameter seems to lie outside (or at the boundary) of most established methods leaving only the method of moments as last weapon as  demonstrated in this paper.

The approach in this paper (based on an asymptotic transfer of the underlying recurrence and computation of all moments as in \cite{Hwang}) has some generality and, together with the analytic tools used by Jacquet in \cite{Ja}, can be applied to other additive shape parameters of G-tries as well in order to obtain similar central limit theorem phenomena. Moreover, it can also be applied to some non-additive shape parameters such as the insertion cost which was also discussed by Jacquet in \cite{Ja}. However, this example is less interesting because mean and variance have the same order and thus, as in the uniform case from Section~\ref{unif}, the contraction method can easily be applied to prove a central limit theorem (this was already pointed out in \cite{Ja}).

In this paper we considered G-tries where the underlying graph G is an $M$-ary tree. Thus, an interesting question is how about the size (and other shape parameters) for G-tries with a general (acyclic) graph G? For this situation, Jacquet in \cite{Ja} computed the mean, however, even the computation of the variance seems to be a considerable challenge; see the recent paper of Jacquet and Magner \cite{JaMa}.

\end{document}